# Studying Pensions Funds Through an Infinite Servers Nodes Network: A Theoretical Problem

M. A. M. Ferreira

M. Andrade

J. A. Filipe

Instituto Universitário de Lisboa (ISCTE-IUL), BRU-IUL, Lisboa, Portugal

manuel.ferreira@iscte.pt

**Abstract**. This study intends to present a representation of a pensions fund through a stochastic network with two infinite servers nodes. With this representation it is allowed to deduce an equilibrium condition of the system with basis on the identity of the random rates expected values, for which the contributions arrive to the fund and the pensions are paid by the fund. In our study a stochastic network is constructed where traffic is represented. This network allows to study the equilibrium in the system and it is admissible to get a balance to a pensions fund. A specific case is studied. When the arrivals from outside at nodes $A$ and $B$ are according to a Poisson process, with rates $\lambda_A$ and $\lambda_B$, respectively, the system may be seen as a two nodes network where the first node is a $M/G/\infty$ queue and second a $M_t/G/\infty$ queue. For this case in the long term the conditions of equilibrium are as follows: $m_A \lambda_A \alpha_A = m_B(\rho \lambda_A + \lambda_B)\alpha_B$. In this formula it is established a relationship among the two nodes. Several examples are given in the study.

## 1. Introduction

Consider two nodes, service centres, $A$ and $B$ both with infinite servers. The traffic through arches $a$ to $e$ is as it is schematized in Figure 1.

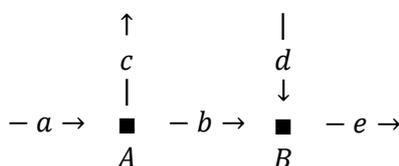

**Figure 1.** Traffic in the stochastic network

The users arrive to node $A$ by arch $a$ at rate $\lambda_A$. And the service time at this node is a positive random variable with distribution function (d.f.) $G_A$ and finite mean $\alpha_A$. After node $A$ the users go to node $B$ through $b$ with probability $p$. Or just abandon the system through arch $c$ with probability $1 - p$.

The users coming directly from outside through $d$ at rate $\lambda_B$ have also access to the service supplied at $B$, according to a positive random variable with d.f. $G_B$ and finite mean $\alpha_B$. The system is abandoned by these users through arch $e$.

In Carvalho [1] this system is suggested as a representation of a pensions fund. So at node *A* arrive individuals that pay, during the service time, their contributions to the fund. The pensioners are at node *B*, which service represents their pensions payment by the fund. This representation reflects also the functions of the common social security funds and that is why it accepts the access of pensioners that have not formerly participated, at node *A*, in the building of the fund.

The target of this study is, having this representation in mind, to obtain results about the transient behavior of the system from the point of view of its equilibrium and autonomy.

## 2. The Fund Equilibrium

Let $N_A(t)$ and $N_B(t)$ be the random variables (r.v.) that represent the number of individuals by time *t* at nodes *A* and *B*, respectively. Consider also the sets of r.v., i.i.d.:

$$X_{A_1}(t), X_{A_2}(t), X_{A_3}(t), \ldots, (X_{B_1}(t), X_{B_2}(t), X_{B_3}(t), \ldots)$$

which designate the unitary contributions, pensions by time *t*, with mean $m_A(t)$ and $m_B(t)$.

The system is in equilibrium when the expected values of the rates at which the contributions are being received and the pensions are being paid by the fund are identical:

$$E\left[\sum_{i=1}^{N_A(t)} X_{A_i}(t)\right] = E\left[\sum_{j=1}^{N_B(t)} X_{B_j}(t)\right].$$

That is, by Wald's equation:

$$m_A(t)E[N_A(t)] = m_B(t)E[N_B(t)] \qquad (1).$$

Eq. (1) just stays that at each instant the mean value of the unitary pension should be proportional to the mean value of the unitary contribution, with the ratio between the averages of the numbers of contributors and pensioners as proportionality factor. Being $t = 0$ the origin time, its solution corresponds, for $t > 0$, to the following pairs:

$$(m_A(t); m_B(t)) = \left(m_A(t); \frac{m_A(t)E[N_A(t)]}{E[N_B(t)]}\right),$$

where $m_A(t)$ is independent of the equilibrium.

If the mean value of the unitary pension is initially 1, and grows continuously with an interest rate *r*,

$$m_B(t) = e^{rt}$$
$$m_A(t) = e^{rt}(E[N_B(t)]/E[N_A(t)]).$$

It is elementary, after Eq. (1),

$$E[N_A(t)] < E[N_B(t)] \Rightarrow m_A(t) > m_B(t).$$

So, in equilibrium, the mean value of the unitary pension is smaller than the mean value of unitary contribution whenever the number of pensioners at *B* is bigger than the number of contributors at *A*.

## 3. The Case of Poisson Arrivals

If the arrivals from outside at nodes A and B are according to a Poisson process, with rates $\lambda_A$ and $\lambda_B$, respectively, the system may be seen as a two nodes network where the first node is a $M/G/\infty$ queue and second a $M_t/G/\infty$ queue, see for instance Ferreira and Andrade [2]. So, $N_A(t)$ is Poisson distributed with parameter, see Ross [3]

$$\lambda_A \int_0^t (1 - G_A(v)) \, dv.$$

The output of the first node is a non-homogeneous Poisson process with intensity function $\lambda_A G_A(t)$ and, consequently, the global arrivals rate at node B is $p\lambda_A G_A(t) + \lambda_B$. Under this conditions $N_B(t)$ is Poisson distributed with parameter, see Harrison and Lemoine [4]:

$$\int_0^t (p\lambda_A G_A(v) + \lambda_B)(1 - G_B(t - v)) dv.$$

And Eq. (1) is written like this

$$m_A(t)\lambda_A \int_0^t (1 - G_A(v)) \, dv = m_B(t) \int_0^t (p\lambda_A G_A(v) + \lambda_B)(1 - G_B(t - v)) dv \quad (2).$$

When $t \to \infty$ the equilibrium conditions assumes the following form where $m_i = \lim_{t\to\infty} m_i(t), i = A, B$:

$$m_A \lambda_A \alpha_A = m_B(\rho\lambda_A + \lambda_B)\alpha_B \quad (3).$$

If the service times at nodes A and B have d.f. concentrated in the intervals $[0, a]$ and $[0, b]$, $m_A\lambda_A\alpha_A = m_B(\rho\lambda_A + \lambda_B)\alpha_B$ for $t \geq a + b$.

## 4. Example
In this section some concrete examples of service times distributions will be considered.

*4.1. Uniformly Distributed Service Times*
If the service times are uniformly distributed, supposing that $\alpha_B < \alpha_A$, it is obtained for Eq. (2) in $0 \leq t < 2\alpha_A + 2\alpha_B$, not to repeat what has just been mentioned:

i) $m_A(t)\lambda_A \left(t - \frac{t^2}{4\alpha_A}\right) = m_B(t)\lambda_B \left(t - \frac{t^2}{4\alpha_B}\right) + m_B(t)p\lambda_A \left(\frac{t^2}{4\alpha_A} - \frac{t^3}{24\alpha_A\alpha_B}\right)$, if $0 \leq \frac{t}{2} < \alpha_B$

ii) $m_A(t)\lambda_A \left(t - \frac{t^2}{4\alpha_A}\right) = m_B(t)\lambda_B\alpha_B + m_B(t)p\lambda_A \left(-\frac{\alpha_B^2}{3\alpha_A} - \frac{t\alpha_B}{2\alpha_A}\right), \alpha_B \leq \frac{t}{2} < \alpha_A$

iii) $m_A(t)\lambda_A\alpha_A = m_B(t)\lambda_B\alpha_B + m_B(t)p\lambda_A \left(-\alpha_A - \frac{\alpha_A^2}{12\alpha_B} + t - \frac{(t-\alpha_A)^2}{4\alpha_B} + \frac{(t-2\alpha_B)^3}{24\alpha_A\alpha_B}\right)$, if $\alpha_A \leq \frac{t}{2} < \alpha_A + \alpha_B$.

*4.2. Exponentially Distributed Service Times*
If the service times are exponentially distributed the equilibrium distribution is given by:

i) $m_A(t)\lambda_A\alpha_A\left(1-e^{-\frac{t}{\alpha_A}}\right) = m_B(t)(p\lambda_A + \lambda_B)\alpha_B\left(1-e^{-\frac{t}{\alpha_B}}\right) - m_B(t)\frac{p\lambda_A\alpha_A\alpha_B}{\alpha_A-\alpha_B}\left(e^{-\frac{t}{\alpha_A}} - e^{-\frac{t}{\alpha_B}}\right)$, if $\alpha_A \neq \alpha_B$

ii) $m_A(t)\lambda_A\alpha_A\left(1-e^{-\frac{t}{\alpha_A}}\right) = m_B(t)(p\lambda_A + \lambda_B)\alpha_A\left(1-e^{-\frac{t}{\alpha_A}}\right) - m_B(t)p\lambda_A t e^{-\frac{t}{\alpha_A}}$, if $\alpha_A = \alpha_B$

*4.3. Service Times with a Particular Distribution Function*

Solving Eq. (2) in the way presented above becomes quite difficult with other standard distributions for the service times. So now it will be considered a collection of d.f.'s, see Ferreira and Andrade [5] and Ferreira and Andrade [6], for the service times given by

$$G_i(v) = 1 - \frac{(1-e^{-\rho_i})(\gamma_i + \beta_i)}{\gamma_i e^{-\rho_i}(e^{(\gamma_i+\beta_i)v} - 1) + \gamma_i}, v \geq 0, \gamma_i > 0, \rho_i > 0, -\gamma_i \leq \beta_i \leq \frac{\gamma_i}{e^{-\rho_i} - 1}, i = A, B.$$

The mean distribution is $\alpha_i = \rho_i/\gamma_i$. In this case Eq. (2) becomes

$$m_A(t)\frac{\lambda_A}{\gamma_A} \ln\frac{e^{(\gamma_A+\beta_A)t}}{e^{-\rho_A}(e^{(\gamma_A+\beta_A)t} - 1) + 1}$$
$$= m_B(t)\frac{p\lambda_A + \lambda_B}{\gamma_B} \ln\frac{e^{(\gamma_B+\beta_B)t}}{e^{-\rho_B}(e^{(\gamma_B+\beta_B)t} - 1) + 1} - m_B(t)p\lambda_A I(t)$$

where

$$I(t) = \int_0^t \frac{(1-e^{-\rho_A})(\gamma_A + \beta_A)}{\gamma_A e^{-\rho_A}(e^{(\gamma_A+\beta_A)v} - 1) + \gamma_A} \times \frac{(1-e^{-\rho_B})(\gamma_B + \beta_B)}{\gamma_B e^{-\rho_B}(e^{(\gamma_B+\beta_B)(t-v)} - 1) + \gamma_B} dv.$$

$I(t)$ is non-negative and not bigger than
$$\frac{(\gamma_A + \beta_A)(\gamma_A + \beta_A)t}{\gamma_A + \gamma_B}.$$

*4.4. Approximations*

The Eq. (2) solution seems to be significantly more complex in circumstances different from those that have been mentioned. For instance, if the service times follow a LogNornal, Gama or Weibull distributions. In some cases, only the numerical solution can eventually be stained.

For appropriate values of *t*, the following approximations concerning the equilibrium conditions are suggested:

$$\frac{m_B(t)}{m_A(t)} \cong \frac{\lambda_A\alpha_A}{(p\lambda_A + \lambda_B)\alpha_B} \quad (4);$$

$$\frac{m_B(t)}{m_A(t)} \cong \frac{\lambda_A}{\lambda_B} \quad (5).$$

Eq. (4) seems reasonable for values of *t* big enough and Eq. (5) is preferred for *t* close to zero. For details see Figueira and Ferreira [7].

## 5. Observations

Some values of the parameters $p$ and $\lambda_B$ have a special influence in the system behaviour. One may consider the suppression of the arch *b* when $p = 0$, of the arch *c* when $p = 1$ or of the arch *d* for $\lambda_B = 0$. Under those circumstances the traffic in those arches can be neglected.

It may be admitted that the ratio $m_B(t)/m_A(t)$ remains constant. This corresponds to the assumption that all the users of the system face identical conditions of effort and benefit, independently of the moment they join the system. Eq. (3) supplies a natural candidate for the value of that constant: $\lambda_A \alpha_A / (p\lambda_A + \lambda_B)\alpha_B$. In such situation Eq. (2) should include an "excess" functions $h(t)$:

$$h(t) = m_B(t) \frac{\lambda_A \alpha_A}{(p\lambda_A + \lambda_B)\alpha_B} \int_0^t (p\lambda_A G_A(v) + \lambda_B)(1 - G_B(t-v)) dv - m_A(t)\lambda_A \int_0^t (1 - G_A(v)) dv.$$

The function $h(t)$ is also interpreted in the sense of the expected value of a random variable depending on *t*. This approach can be generalized in a natural way to some other predefined function $m_B(t)/m_A(t)$.

Assuming that the system is initially empty appears to be a strong restriction of the analysis performed. When someone meets the system already in operation and does not known when it did start, the results that have been mentioned seem to have a lesser utility. In such case, there re-evaluation or finding a estimation procedure for the initial time are determinant for practical purposes.

**Acknowledgments**

This work was financially supported by FCT through the Strategic Project PEst-OE/EGE/UI0315/2011.

**References**
[1]    Carvalho, P., "*Planos e Fundos de Pensões*", Texto Editora, Lisboa, (1993).
[2]    Ferreira, M. A. M. and Andrade, M., "*Fundaments of Theory of Queues*", International Journal of Academic Research, Vol 3, No 1, Part II, pp. 427-429, (2011).
[3]    Ross, S., *Stochastic Processes,* 2nd Edition, Wiley, New York, (1996).
[4]    Harrison, J. and Lemoine, A., *"A Note on Networks of Infinite-Server Queues"*, Journal of Applied Probability, Vol 18, pp. 561-567, (1981).
[5]    Ferreira, M. A. M. and Andrade, M., "*The Ties Between the $M/G/\infty$ Queue System Transient Behaviour and the Busy Period*", International Journal of Academic Research, Vol 1, No 1, pp. 84-92, (2009).
[6]    Ferreira, M. A. M. and Andrade, M., "*Looking to a $M/G/\infty$ System Occupation Through a Ricatti Equation*", Journal of Mathematics and Technology, Vol 1, No 2, pp. 58-62, (2010).
[7]    Figueira, J. and Ferreira, M. A. M., *"Representation of a Pensions Fund by a Stochastic Network with Two Nodes: an Exercise"*, Portuguese Review of Financial Markets, Vol 2, No 1, pp. 75-81, (1999).
[8]    Andrade, M. "*A Note on Foundations of Probability*", Journal of Mathematics and Technology, Vol 1, No 1, pp. 96-98, (2010).

[9]     Andrade, M., Ferreira, M. A. M., Filipe, J. A. And Coelho, M., *"The Study of a Pensions Fund Equilibrium Through an Infinite Servers Nodes Network"*, International Journal of Academic Research, Vol 4, No 3, pp. 205-208, (2012).
[10]   Ferreira, M. A. M., *"A Note on Jackson Networks Sojourn Times"*, Journal of Mathematics and Technology, Vol 1, No 1, pp. 91-95, (2010).
[11]   Ferreira, M. A. M. and Andrade, M., *"An Infinite Servers Nodes Network in the Study of a Pensions Fund"*, International Journal of Latest Trends in Finance and Economic Sciences, Vol. 1 (2), pp 88-90, (2011).
[12]   Figueira, J., *"Aplicação dos Processos de Difusão e da Teoria do Renovamento num Estudo de Reservas Aleatórias"*, PhD Thesis, ISCTE, (2003).